\documentclass{amsart}
\usepackage{amsmath}
\usepackage{amssymb}
\usepackage{amsthm}
\usepackage[all]{xy}

\newtheorem{theorem}{Theorem}[section]
\newtheorem{corollary}[theorem]{Corollary}

\newtheorem{lemma}[theorem]{Lemma}
\newtheorem{proposition}[theorem]{Proposition}

\def\Q{\mathbb{Q}}
\def\Z{\mathbb{Z}}
\def\op{\operatorname}

\def\aut{\op{Aut}}

\def\da{\op{\delta}}
\def\dm{\op{dim}}
\def\endo{\op{End}}
\def\eps{\epsilon}
\def\f{\mathcal{F}}
\def\fm{f^{-1}}
\def\fs{f^{\sg}}
\def\fp{\mathbb F_p}
\def\fq{\mathbb F_q}
\def\ft{\mathbb F_2}
\def\gal{\op{Gal}}
\def\imp{\,\Longrightarrow\,}
\def\kb{\overline{k}}
\def\la{\lambda}
\def\lra{\longrightarrow}
\def\md#1{\ \mbox{\rm(mod }{#1})}
\def\pr#1{\mathbb P^{#1}}
\def\rk{\op{rk}}
\def\sg{\sigma}
\def\sgn#1{\renewcommand\arraystretch{.8}\left(\!\begin{tabular}{c}$\!\!-1$\\\hline $#1$\end{tabular}\!\right)\renewcommand\arraystretch{1.}}

\def\tq{\,\,|\,\,}
\def\tr{\op{tr}}

\begin{document}

\title{Zeta Function and Cryptographic Exponent of Supersingular Curves of Genus 2}

\author{Gabriel Cardona}

\address{Dept. Ci\`encies Matem\`atiques i Inform\`atica,
  Universitat de les Illes Balears,
  07122, Palma de Mallorca, Spain}
  
  \email{gabriel.cardona@uib.es}

\author{Enric Nart}

\address{Departament de Matem\`atiques,
        Universitat Aut\`onoma de Barcelona,
Edifici C,
        08193 Bellaterra, Barcelona, Spain}
\email{nart@mat.uab.cat}

\thanks{The authors acknowledge support from the projects MTM2006-15038-C02-01 and  MTM2006-11391 from the Spanish MEC}

\begin{abstract}
We compute in a direct (not algorithmic) way the zeta function of all supersingular curves of genus $2$ over a finite field
$k$, with many geometric automorphisms. We display these computations in an appendix where we select a family of representatives of all these curves up to $\kb$-isomorphism and we exhibit equations and the zeta function of all their $\kb/k$-twists. As an application we obtain a direct computation of the cryptographic exponent of the Jacobians of these curves.       
\end{abstract}

\maketitle

\section*{Introduction}
One-round tripartite Diffie-Hellman, identity based encryption, and short digital signatures are some problems for which good solutions have recently been found, making critical use of pairings on supersingular abelian varieties over a finite field $k$. The cryptographic exponent $c_A$ of a supersingular abelian variety $A$ is a half-integer that measures the security against an attack on the DL problem based on the Weil or the Tate pairings. Also, it is relevant to determine when pairings can be efficiently computed. Rubin and Silverberg showed in \cite{rubsil} that this invariant is determined by the zeta function of $A$.  

In this paper we give a direct, non-algorithmic procedure to compute the zeta function of a supersingular curve of genus $2$, providing thus a direct computation of the cryptographic exponent of its Jacobian. This is achieved in Sect. \ref{secu}. For even characteristic the results are based on \cite{mn2} and are summarized in Table \ref{t2}; for odd characteristic we use results of Xing and Zhu on the structure of the group of $k$-rational points of a supersingular abelian surface and we almost determine the zeta function in terms of the Galois structure of the set of Weierstrass points of the curve (Tables \ref{t3}, \ref{t4}). In the rest of the paper we obtain a complete answer in the case of curves with many automorphisms. In Sect. \ref{sec2} we study the extra information provided by these automorphisms and we show how to obtain the relevant data to compute the zeta funtion of a twisted curve in terms of data of the original curve and the $1$-cocycle defining the twist. In Sect. \ref{sec3} we select a family of representatives of these curves up to $\kb$-isomorphism and we apply the techniques of the previous section to deal with each curve and all its $\kb/k$-twists. The results are displayed in an Appendix in the form of tables.   

In what cryptographic applications of pairings concerns, curves with many automorphisms are interesting too because they are natural candidates to provide distortion maps on the Jacobian. In this regard the computation of the zeta function is a necessary step to study the structure of the endomorphism ring of the Jacobian (cf. \cite{GPRS}).  

\section{Zeta Function and Cryptographic Exponent}\label{secu}
Let $p$ be a prime number and let $k=\fq$ be a finite field of characteristic $p$. We denote by $k_n$ the extension of degree $n$ of $k$ in a fixed algebraic closure $\kb$, $G_k:=\gal(\kb/k)$ is the absolute Galois group of $k$, and
$\sg\in G_k$ the Frobenius automorphism. 

Let $C$ be a projective, smooth, geometrically irreducible, supersingular curve of genus $2$ defined over $k$. The Jacobian $J$ of $C$ is a supersingular abelian surface over $k$ (the $p$-torsion subgroup of $J(\kb)$ is trivial).
Let us recall how supersingularity is reflected in a model of the curve $C$:

\begin{theorem}\label{sscrit}
If $p$ is odd, any curve of genus $2$ defined over $k$ admits an affine Weierstrass model $y^2=f(x)$, with $f(x)$ a separable polynomial in $k[x]$ of degree $5$ or $6$. The curve is supersingular if and only if $M^{(p)}M=0$, where $M$, $M^{(p)}$ are the matrices:
$$
M=\begin{pmatrix}c_{p-1}&c_{p-2}\\c_{2p-1}&c_{2p-2}\end{pmatrix},\quad
M^{(p)}=\begin{pmatrix}c_{p-1}^p&c_{p-2}^p\\c_{2p-1}^p&c_{2p-2}^p\end{pmatrix},\quad f(x)^{(p-1)/2}=\sum_{j\ge0}c_jx^j\enspace.
$$  

If $p=2$ a curve of genus $2$ defined over $k$ is supersingular if and only if it admits an affine Artin-Schreier model  $y^2+y=f(x)$, with $f(x)$ an arbitrary polynomial in $k[x]$ of degree $5$.
\end{theorem}

For the first statement see \cite{yui} or \cite{iko}, for the second see \cite{vdG}.  

For any simple supersingular abelian variety $A$ defined over $k$, Rubin and Silverberg computed in \cite{rubsil} the {\it cryptographic exponent}  $c_A$, defined as the half-integer such that $q^{c_A}$ is the size of the smallest field $F$ such that every cyclic subgroup of $A(k)$ can be embedded in $F^*$. This invariant refines the concept of {\it embedding degree}, formerly introduced as a measure of the security of the abelian variety against the attacks to the DLP by using the Weil pairing \cite{mov} or the Tate pairing \cite{fruc} (see for instance \cite{gal}).
   
Let us recall the result of Rubin-Silverberg, adapted to the dimension two case. After classical results of Tate and Honda, the isogeny class of $A$ is determined by the Weil polynomial of $A$, $f_A(x)=x^4+rx^3+sx^2+qrx+q^2\in\Z[x]$, which is the characteristic polynomial of the Frobenius endomorphism of the surface. For $A$ supersingular the roots of $f_A(x)$ in $\overline{\Q}$ are of the form $\sqrt{q}\,\zeta$, where $\sqrt{q}$ is the positive square root of $q$ and $\zeta$ is a primitive $m$-th root of unity. 

\begin{theorem}
Suppose $A$ is a simple supersingular abelian surface over $\fq$ and let $\ell>5$ be any prime number dividing $|A(\fq)|$.
Then, the smallest half-integer $c_A$ such that $q^{c_A}-1$ is an integer divisible by $\ell$ is given by
$$
c_A=\left\{\begin{array}{ll}m/2,&\qquad \mbox{if $q$ is a square},\\m/(2,m),&\qquad \mbox{if $q$ is not a square}\enspace.
\end{array}\right.
$$
\end{theorem}

In particular, the cryptographic exponent $c_A$ is an invariant of the isogeny class of $A$.  The complete list of simple supersingular isogeny classes of abelian surfaces can be found in \cite[Thm. 2.9]{mn}. It is straightforward to find out the $m$-th root of unity in each case. We display the computation of $c_A$ in Table \ref{t1}. 
\begin{table}
\caption{Cryptographic exponent $c_A$ of the simple supersingular
abelian surface $A$ with Weil polynomial $f_A(x)=x^4+rx^3+sx^2+qrx+q^2$}\label{t1}
\begin{center}
\begin{tabular}{|c|l|c|}
\hline $(r,s)$&conditions on $p$ and $q$ &$c_A$\\ \hline 
 $(0,-2q)$& $q$ nonsquare
&$1$\\$(0,2q)$& $q$ square, $p\equiv 1\md4$ &$2$\\ $(2\sqrt q,3q)$& $q$ square, $p\equiv 1\md3$&$3/2$\\ $(-2\sqrt q,3q)$& $q$ square, $p\equiv 1\md3$&$3$
\\ 
 $(0,0)$& ($q$ nonsquare, $p\ne 2$) or
($q$ square, $p\not\equiv 1\md8$)&$4$\\ 
 $(0,q)$& $q$ nonsquare& $3$\\ $(0,-q)$& ($q$ nonsquare, $p\ne 3$) or ($q$ square, $p\not\equiv
1\md{12}$)& $6$\\ $(\sqrt q,q)$& $q$ square, $p\not\equiv 1\md5$&$5/2$\\ $(-\sqrt q,q)$& $q$ square, $p\not\equiv 1\md5$&$5$\\ 
$(\pm\sqrt{5q},3q)$& $q$ nonsquare, $p=5$ &$5$\\ $(\pm\sqrt
{2q},q)$& $q$ nonsquare, $p=2$&$12$\\
\hline
\end{tabular}
\end{center}
\end{table}

Therefore, the computation of the cryptographic exponent of the Jacobian $J$ of a supersingular curve $C$ amounts to the computation of the Weil polynomial of $J$, which is related in a well-known way to the zeta function of $C$. We shall call $f_J(x)$ the {\it Weil polynomial of $C$} too.

The computation of $f_J(x)$ has deserved a lot of attention because for the cryptographic applications one needs to know the cardinality $|J(\fq)|=f_J(1)$ of the group of rational points of the Jacobian. However, in the supersingular case the current ``counting points" algorithms are not necessary because there are more direct ways to compute the polynomial $f_J(x)$. 

The aim of this section is to present these explicit methods, which take a different form for $p$ odd or even. For $p=2$ the computation of $f_J(x)$ is an immediate consequence of the methods of \cite{mn2}, based on ideas of van der Geer-van der Vlugt; for $p>2$ we derive our results from the group structure of $J(\fq)$, determined in \cite{xin}, \cite{zhu}, and from the exact knowledge of what isogeny classes of abelian surfaces do contain Jacobians \cite{hnr}. In both cases we shall show that $f_J(x)$ is almost determined by the structure as a Galois set of a finite subset of $\kb$, easy to compute from the defining equation of $C$.

\subsection{Computation of the Zeta Function when $p=2$}\label{secuu}
We denote simply by $\tr$ the absolute trace
$\tr_{k/\ft}$. Recall that $\ker(\tr)=\{x+x^2\tq x\in k\}$ is an $\ft$-linear subspace of $k$ of codimension $1$.

Every projective smooth geometrically irreducible supersingular curve $C$ of genus $2$ defined over $k$ admits an affine Artin-Schreier
model of the type:
$$
C\colon\quad y^2+y=ax^5+bx^3+cx+d,\quad a\in k^*,\,b,\,c,\,d\in
k,
$$
which has only one point at infinity  \cite{vdG}. The change of variables $y=y+u$, $u\in k$, allows us to suppose that 
 $d=0$ or $d=d_0$, with $d_0\in
k\setminus\ker(\tr)$ fixed. Twisting $C$ by the hyperelliptic twist consists in adding $d_0$ to the defining equation. If we
denote by $J'$ the Jacobian of the twisted curve we have $f_{J'}(x)=f_J(-x)$. Thus, for the computation of $f_J(x)$ we can assume that $d=0$.

The structure as a $G_k$-set of the set of roots in $\kb$ of the polynomial $P(x)=a^2x^5+b^2x+a\in k[x]$ almost determines the zeta function of $C$ \cite[Sect.3]{mn2}. 
\begin{table}
\caption{Weil polynomial $x^4+rx^3+sx^2+qrx+q^2$ of the curve $y^2+y=ax^5+bx^3+cx$, for $q$ nonsquare (left) and $q$ square (right)}
\label{t2}
\begin{center}
\begin{tabular}{|c|c|c|}
\hline $P(x)$&$N,\,M$ &$(r,s)$\\ \hline 
 $(1)(4)$& $N=0$
&$(\pm\sqrt{2q},2q)$\\& $N=1$ &$(0,0)$\\ 
 \hline 
 $(2)(3)$& $M=0$
&$(\pm\sqrt{2q},q)$\\& $M=1$ &$(0,q)$\\ 
\hline
 & $N=0$
&$(\pm2\sqrt{2q},4q)$\\$(1)^3(2)$& $N=1$ &$(0,2q)$\\ 
 & $N=2$
&$(0,0)$\\& $N=3$ &$(0,-2q)$\\ 
\hline
\end{tabular}\qquad\qquad
\begin{tabular}{|c|c|c|}
\hline $P(x)$&$N,\,M$ &$(r,s)$\\ \hline 
 $(5)$& 
&$(\pm\sqrt{q},q)$\\\hline & $N=0$ &$(0,-q)$\\
$(1)^2(3)$& $N=1$ &$(0,q)$\\& $N=2$ &$(\pm 2\sqrt{q},3q)$\\
 \hline 
 & $M=0$
&$(\pm 2\sqrt{q},2q)$\\ 
 $(1)(2)^2$& $M=1$
&$(0,0)$\\ 
 & $M=2$
&$(0,2q)$\\ 
\hline
 & $N=1$
&$(0,-2q)$\\ 
 $(1)^5$& $N=3$
&$(0,2q)$\\ 
 & $N=5$
&$(\pm4\sqrt{q},6q)$\\ 
\hline
\end{tabular}
\end{center}
\end{table}

In Table \ref{t2} we write $P(x)=(n_1)^{r_1}(n_2)^{r_2}\cdots(n_m)^{r_m}$
to indicate that $r_i$ of the irreducible factors of $P(x)$ have degree $n_i$. Also, we consider the linear operator $T(x):=\tr((c+b^2a^{-1})x)$ and we define
$$
\begin{array}{l}
N:=\mbox{number of roots $z\in k$ of $P(x)$ s.t. }T(z)=0,\\
M:=\mbox{number of irred. quadratic factors $x^2+vx+w$ of $P(x)$ s.t. }T(v)=0\enspace.
\end{array}
$$

The ambiguity of the sign of $r$ can be solved by computing $nD$ in the Jacobian, where $n$ is one of the presumed values of $|J(\fq)|$ and $D$ is a random rational divisor of degree $0$.     

\subsection{Computation of the Zeta Function when $p$ is odd}
Let $A$ be a supersingular abelian surface over $k$ and let $\rk_2(A):=\dm_{\ft}(A[2](k))$.

The structure of $A(k)$ as an abelian group was studied in \cite{xin}, \cite{zhu}, where it is proven that it is almost determined by the isogeny class of $A$. In fact, if $F_i(x)$ are the different irreducible factors of $f_A(x)$ in $\Z[x]$:
$$
f_A(x)=\prod_{i=1}^sF_i(x)^{e_i},\ 1\le s \le 2\quad \imp\quad 
A(k)\simeq \oplus_{i=1}^s\left(\Z/F_i(1)\Z\right)^{e_i},
$$
except for the following cases:

(a) $p\equiv 3\md4$, $q$ is not a square and $f_A(x)=(x^2+q)^2$,

(b) $p\equiv 1\md4$, $q$ is not a square and $f_A(x)=(x^2-q)^2$.

(c) $q$ is a square and $f_A(x)=(x^2-q)^2$.

The possible structure of $A(k)$ in cases  (a) and (b)
is:
$$
A(k)\simeq \left(\Z/F(1)\Z\right)^m\oplus \left(\Z/(F(1)/2)\Z\oplus\Z/2\Z\right)^n,
$$
where $F(x)$ denotes respectively $x^2+q$, $x^2-q$, and $m,\,n$ are non-negative integers such that $m+n=2$ \cite[Thm. 1.1]{zhu}.
In case (c) we have either:
$$
\begin{array}{l}
A(k)\simeq \left(\Z/((q-1)/2)\Z\right)^2\oplus \left(\Z/2\Z\right)^2,\qquad \mbox{ or }\\A(k)\simeq\left(\Z/((q-1)/2^m)\Z\right)\oplus\left(\Z/((q-1)/2^n)\Z\right)\oplus\left(\Z/2^{m+n}\Z\right),
\end{array}
$$where $0\le m,\,n\le v_2(q-1)$ \cite[Thm. 3]{xin}. In this last case we have $\rk_2(A)>1$; in fact, $v_2(1-\sqrt{q})+v_2(1+\sqrt{q})=v_2(1-q)=(1/2)v_2(F(1))$ and we can apply \cite[Lem. 4]{xin}
to conclude that $A(k)$ has a subgroup isomorphic to $(\Z/2\Z)^2$.

Consider now a supersingular curve $C$ of genus 2 defined over $k$, given by a Weierstrass equation $y^2=f(x)$, for some separable polynomial $f(x)\in k[x]$ of degree $5$ or $6$. Let $J$ be its Jacobian variety, $W=\{P_0,P_1,P_2,P_3,P_4,P_5\}\subseteq C(\kb)$ the set of Weierstrass points of $C$, and $W(k)\subseteq W$ the subset of
$k$-rational Weierstrass points. Our aim is to show that the structure of $W$ as a $G_k$-set contains enough information on the $2$-adic value of $|C(k)|$ and $|J(k)|$ to almost determine the 
 polynomial $f_J(x)=x^4+rx^3+sx^2+qrx+q^2$.

From the fundamental identities
$$
|C(k)|=q+1+r,\qquad |J(k)|=f_J(1)=(q^2+1)+(q+1)r+s,
 $$ and the free action of the hyperelliptic involution on $C(k)\setminus W(k)$ we get
\begin{equation}\label{mod2}
r\equiv |W(k)|\md2,\qquad s\equiv |J(k)|\md2\enspace.
\end{equation}
On the other hand, $J[2]$ is represented by the classes of the 15 divisors:
$$P_i-P_0, \ \ 1\le i\le 5,\qquad  \mbox{and }\ \quad P_i+P_j-2P_0, \ \ 1\le i< j\le 5,
$$together with the trivial class. 

\begin{lemma}
Let $D=P_i-P_j$, with $i\ne j$, or $D=P_i+P_j-2P_0$, with $0,\,i,\,j$ pairwise different. Then, the class of the divisor $D$ is $k$-rational if and only if $P_i,\,P_j$ are both $k$-rational or quadratic conjugate. 
\end{lemma}

Hence, the Galois structure of $W$ determines $\rk_2(J)$ and this limits the possible values of the zeta function of $C$. Our final results are given in Tables \ref{t3},\,\ref{t4}, where we write $W= (n_1)^{r_1}(n_2)^{r_2}\cdots(n_m)^{r_m}$ to indicate that there are $r_i$ $G_k$-orbits of length $n_i$ of Weierstrass points. If $f(x)$ has degree $6$ this Galois structure mimics the decomposition $f(x)=(n_1)^{r_1}(n_2)^{r_2}\cdots(n_m)^{r_m}$ (same notation as in Sect. \ref{secuu}) of $f(x)$ into a product of irreducible polynomials $k[x]$. If $f(x)$ has degree $5$ then $W=(1)f(x)$, because in these models the
point at infinity is a $k$-rational Weierstrass point. 
\begin{table}
\caption{Weil polynomial $x^4+rx^3+sx^2+qrx+q^2$ of the curve $C$ when $q$ is nonsquare. The sign $\eps$ is the Legendre symbol $(-1/p)$ }\label{t3}
\begin{center}
\begin{tabular}{|c|c|c|c|}
\hline $W$&  $p$&$\rk_2(J)$&$(r,s)$\\
\hline$(1)^6\mbox{ or }(1)^4(2)$&&$4,\,3$&$(0,-2\eps q)$\\
\hline$(1)^2(2)^2$\mbox{ or }$(2)^3$&&$2$&$(0,\pm2q)$\\
\hline$(1)^3(3)$&&$2$&\mbox{\scriptsize not possible}\\\hline
$(1)(2)(3)$&$p>3$&$1$&\mbox{\scriptsize not possible}\\
&$p=3$&&$(\pm\sqrt{3q},2q)$\\\hline
$(1)^2(4)$\mbox{ or }$(2)(4)$&&$1$&$(0,0)$\\\hline
$(1)(5)$&$p\ne5$&$0$&\mbox{\scriptsize not possible}\\
&$p=5$&&$(\pm\sqrt{5q},3q)$\\\hline
&$p\equiv1\md3$&&$(0,q)$\\
$(3)^2$&$p\equiv-1\md3$&$0$&$(0,\eps q)$\\
&$p=3$&&\mbox{\scriptsize not possible}\\\hline
$(6)$&$p\equiv-1\md3$&$0$&$(0,\pm q)$\\
&$p\not\equiv-1\md3$&&$(0,q)$\\\hline
\end{tabular}
\end{center}
\end{table}
\begin{table}
\caption{Weil polynomial $x^4+rx^3+sx^2+qrx+q^2$ of the curve $C$ when $q$ is a square}\label{t4}
\begin{center}
\begin{tabular}{|c|c|c|c|}
\hline $W$&  $p$&$\rk_2(J)$&$(r,s)$\\
\hline$(1)^6$&&$4$&$(0,-2q)\mbox{ or }(\pm4\sqrt{q},6q)$\\
\hline$(1)^4(2)$&&$3$&$(0,-2q)$\\
\hline$(1)^2(2)^2$\mbox{ or }$(2)^3$&&$2$&$(0,\pm2q)$\\\hline
$(1)^3(3)$&$p>3$&$2$&\mbox{\scriptsize not possible}\\
&$p=3$&&$(\pm\sqrt q,0)$\\\hline
$(1)(2)(3)$&&$1$&\mbox{\scriptsize not possible}\\\hline
$(1)^2(4)$\mbox{ or }$(2)(4)$&$p\equiv1\md8$&$1$&\mbox{\scriptsize not possible}\\
&$p\not\equiv1\md8$&&$(0,0)$\\\hline
$(1)(5)$&$p\equiv1\md5$&$0$&\mbox{\scriptsize not possible}\\
&$p\not\equiv1\md5$&&$(\pm\sqrt{q},q)$\\\hline
$(3)^2$&&&$(0,q)$\mbox{\scriptsize  \,or }$(\pm 2\sqrt{q},3q)$\\\hline
$(6)$&$p\equiv5\md{12}$&$0$&$(0,\pm q)$\\
&$p\not\equiv5\md{12}$&&$(0,q)$\\\hline
\end{tabular}
\end{center}
\end{table}

The proof of the content of Tables 3 and 4 is elementary, but long. Instead of giving all details we only sketch the main ideas:

(I) Waterhouse determined all possible isogeny classes of supersingular elliptic curves \cite{Wa}. Thus, it is possible to write down all isogeny classes of supersingular abelian surfaces by adding to the simple classes given in Table \ref{t1} the split isogeny classes. By \cite{hnr} we know exactly what isogeny classes of abelian surfaces do not contain Jacobians and they can be dropped from the list. By the results of Xing and Zhu we can distribute the remaining isogeny classes according to the possible values of $\rk_2$.

(II) Each structure of $W$ as a $G_k$-set determines the value of $\rk_2$ and, after (I), it has a reduced number of possibilities for the isogeny classes. By using (\ref{mod2}) and looking for some incoherence in the behaviour under scalar extension to $k_2$ or $k_3$ of both, the Galois structure of $W$ and the possible associated isogeny classes, we can still discard  some of these possibilities.

In practice, among the few possibilities left in Tables 3 and 4 we can single out the isogeny class of the Jacobian of any given supersingular curve by computing iterates of random divisors of degree zero. However, if $C$ has many automorphisms they provide enough extra information to completely determine the zeta function. This will be carried out in the rest of the paper. In the Appendix we display equations of the supersingular curves with many automorphisms and their Weil polynomial. 

\section{Zeta Function of Twists}\label{sec2}
In this section we review some basic facts about twists and we  show how to compute different properties of a twisted curve in terms of the defining $1$-cocycle. From now on the ground field $k$ will have odd characteristic.

Let $C$ be a supersingular curve of genus $2$ defined over $k$ and let $W\subseteq C(\kb)$ be the set of Weierstrass points of $C$. We denote by $\aut(C)$ the $k$-automorphism group of $C$ and by
$\aut_{\kb}(C)$ the full automorphism group of $C$.  

Let $\phi\colon C\lra \pr1$ be a fixed $k$-morphism of degree $2$ and consider the group of reduced geometric automorphisms of $C$:
$$
\aut'_{\kb}(C):=\{u'\in \aut_{\kb}(\pr1)\tq u'(\phi(W))=\phi(W)\}\enspace. 
$$We denote by $\aut'(C)$ the subgroup of reduced automorphisms defined over $k$.

Any automorphism $u$ of $C$ fits into a commutative
  diagram:
  \begin{center}
\leavevmode
\xymatrix{ C \ar[r]^u \ar[d]^\phi & C \ar[d]^\phi \\
\pr1 \ar[r]^{u'} & \pr1
}
\end{center}
  for certain uniquely determined reduced automorphism $u'$. The map $\ u\mapsto u'$ is a
group  homomorphism (depending on $\phi$) and we have a central exact sequence of
 groups compatible with Galois action:
  $$
  1\lra\{1,\iota\} \lra \aut_{\kb}(C) \stackrel{\phi}\lra\aut'_{\kb}(C)\lra 1,
  $$
where $\iota$ is the hyperelliptic involution. This leads to a long exact sequence of Galois cohomology sets:
\begin{multline}\label{cohom}
1\to\{1,\iota\} \to \aut(C) \stackrel{\phi}\to\aut'(C)\stackrel{\da}\to H^1(G_k,\{1,\iota\})\to  H^1(G_k,\aut_{\kb}(C))\to\\ \to
H^1(G_k,\aut'_{\kb}(C))\to
H^2(G_k,\{1,\iota\})\simeq \op{Br}_2(k)=0\enspace.
\end{multline}

The $\kb/k$-twists of $C$ are parameterized by the pointed set $H^1(G_k,\aut_{\kb}(C))$
and, since $k$ is a finite field, a 1-cocycle is determined just by the choice of an automorphism $v\in \aut_{\kb}(C)$. The twisted curve $C_v$ associated to $v$ is defined over $k$ and is determined, up to $k$-isomorphism, by the existence of a
$\kb$-isomorphism
$f\colon C\lra C_v$, such that $\fm\fs=v$.

For instance, the choice $v=\iota$ corresponds to the hyperelliptic twist $C'$; if $C$ is given by a Weierstrass equation $y^2=f(x)$ then $C'$ admits the model $y^2=tf(x)$, for $t\in k^*\setminus (k^*)^2$. We say that $C$ is {\it self-dual} if it is $k$-isomorphic to its hyperelliptc twist. 
If $f_J(x)$ is the Weil polynomial of $C$, the Weil polynomial of $C'$ is $f_{J'}(x)=f_J(-x)$; in particular, for a self-dual curve one has $f_J(x)=x^4+sx^2+q^2$ for some integer $s$.

It is easy to deduce from (\ref{cohom}) the following criterion for self-duality:

\begin{lemma}\label{prima}
The curve $C$ is self-dual if and only if $|\aut'(C)|=|\aut(C)|$.
\end{lemma}

One can easily compute the data $\aut(C_v)$, $\aut'(C_v)$ of the twisted curve $C_v$, in terms of $\aut_{\kb}(C)$, $\aut'_{\kb}(C)$ and the 1-cocycle $v$. 

Let $f\colon C\lra C_v$ be a geometric isomorphism such that $\fm\fs=v$. We have $\aut_{\kb}(C_v)=f\aut_{\kb}(C)\fm$, and the $k$-automorphism group is
\begin{equation}\label{autprima}
\aut(C_v)=\{fu\fm\tq u\in \aut_{\kb}(C),\, u\,v=v\,u^{\sg}\}\enspace.
\end{equation}
Once we fix any $k$-morphism of degree two, $\phi_v\colon C_v\lra\pr1$, it determines a unique geometric automorphism $f'$ of $\pr1$ such that $\phi_vf=f'\phi$. The reduced group of $k$-automorphisms of $C_v$ is 
\begin{equation}\label{kerdelta}
\aut'(C_v)=\{f'u'(f')^{-1}\tq u'\in \aut'_{\kb}(C),\, u'\,v'=v'(u')^{\sg}\}\enspace.
\end{equation} 

In order to compute the zeta function of $C_v$ we consider the geometric isomorphism $f\colon J\lra J_v$ induced by $f$. We still have $\fm\fs=v_{\ast}$, where $v_{\ast}$ is the automorphism of $J$ induced by $v$. Clearly, $\pi_vf=\fs\pi$, where $\pi,\,\pi_v$ are the respective $q$-power Frobenius endomorphisms of $J,\,J_v$. Hence, 
$$
f^{-1}\pi_vf=f^{-1}\fs(\fs)^{-1}\pi_vf=v_{\ast}\pi\enspace.
$$ 
In particular, $\pi_v$ has the same characteristic polynomial than $v_{\ast}\pi$. From this fact one can deduce two crucial results (cf. \cite[Props.13.1,13.4]{hnr}).

\begin{proposition}\label{qsq}
Suppose $q$ is a square. Let $C$ be a supersingular genus $2$ curve over $k$ with Weil polynomial $(x+\sqrt{q})^4$ and let  $v$ be a geometric automorphism of $C$, $v\ne 1,\,\iota$. Then, the Weil polynomial $x^4+rx^3+sx^2+rx+q^2$ of $C_v$ is determined as follows in terms of $v$
(in the column $v^6=1$ we suppose $v^2\ne1,\,v^3\ne1,\iota$):
\begin{center}
\begin{tabular}{|c|c|c|c|c|c|c|c|c|c|}
\hline $v$& $v^2\!=1$& $v^2\!=\iota$&$v^3\!=1$& $v^3\!=\iota$& $v^4\!=\iota$&$v^5\!=1$& $v^5\!=\iota$& $v^6\!=1$& $v^6\!=\iota$\\
\hline $(r,s)$&$(0,\!-2q)$&$(0,2q)$&$(-2\sqrt{q},3q)$&$(2\sqrt{q},3q)$&$(0,0)$&$(-\sqrt{q},q)$&$(\sqrt{q},q)$&$(0,q)$&$(0,\!-q)$\\
\hline
\end{tabular}
\end{center}
\end{proposition}
\begin{proposition}\label{qnsq}
Suppose $q$ is nonsquare. Let $C$ be a supersingular genus $2$ curve over $k$  with Weil polynomial $(x^2+\eps q)^2$, $\eps\in\{1,-1\}$, and let  $v$ be a geometric automorphism of $C$. Then, the Weil polynomial $x^4+rx^3+sx^2+rx+q^2$ of $C_v$ is determined as follows in terms of the order $n$ of the automorphism $vv^{\sg}$:
\begin{center}
\begin{tabular}{|c|c|c|c|c|c|}
\hline $n$&$1$&$2$&$3$&$4$&$6$\\
\hline $(r,s)$&$(0,2\eps q)$&$(0,-2\eps q)$&$(0,-\eps q)$&$(0,0)$&$(0,\eps q)$\\
\hline
\end{tabular}
\end{center}
\end{proposition}

In applying these results the transitivity property of twists can be helpful.

\begin{lemma}\label{trans}
Let $u,v$ be automorphisms of $C$ and let $f\colon C\to C_v$ be a geometric isomorphism with $\fm\fs=v$. Then the curve $C_u$ is the twist of $C_v$ associated to the automorphism $f uv^{-1}\fm$ of $C_v$.
\end{lemma}

For a curve with a large $k$-automorphism group the following remark, together with Tables \ref{t3} and  \ref{t4}, determines in some cases the zeta function:
\begin{lemma}\label{modauto}
Let $\f\subseteq C(k)$ be the subset of $k$-rational points of $C$ that are fixed by some non-trivial $k$-automorphism of $C$. Then,
$$
|C(k)|\equiv |\f|\md{|\aut(C)|}\enspace.
$$
\end{lemma}

\begin{proof}
The group $\aut(C)$ acts freely on $C(k)\setminus \f$.
\end{proof}

Note that $\f$ contains the set $W(k)$ of $k$-rational Weierstrass points, all of them fixed by the hyperelliptic involution $\iota$ of $C$.

In order to apply this result to the twisted curve $C_v$ we need to compute the $G_k$-set structure of $W_v$ and $|\f_v|$ solely in terms of $v$.

\begin{lemma}\label{orbit}
(1) For any $P\in W$ the length of the $G_k$-orbit of $f(P)\in W_v$ is the minimum positive integer $n$ such that $v\,v^{\sg}\cdots v^{\sg^{n-1}}(P^{\sg^n})=P$. In particular, $|W_v(k)|=|\{P\in W\tq v(P^{\sg})=P\}|$.

(2) The map $\fm$ stablishes a bijection between $\f_v$ and the set 

\noindent$\{P\in C(\kb)\tq v(P^{\sg})=P=u(P)\mbox{ for some }1\ne u\in\aut_{\kb}(C),\mbox{ s.t. } u\,v=v\,u^{\sg}\}$.
\end{lemma}

\section{Supersingular curves with many automorphisms}\label{sec3}
For several cryptographic applications of the Tate pairing the use of distortion maps is essential. A distortion map is an endomorphism $\psi$ of the Jacobian $J$ of $C$ that provides an input for which the value of the pairing is non-trivial: $e_{\ell}(D_1,\psi(D_2))\ne 1$ for some fixed $\ell$-torsion divisors $D_1$, $D_2$. The existence of such a map is guaranteed, but in practice it is hard to find it in an efficient way. Usually, one can start with a nice curve $C$ with many automorphisms, consider a concrete automorphism $u\ne 1$, $u\ne\iota$, and look for a distortion map $\psi$ in the subring $\Z[\pi,u_{\ast}]\subseteq \endo(J)$, where $\pi$ is the Frobenius endomorphism of $J$ and $u_{\ast}$ is the automorphism of the Jacobian induced by $u$. If $\Z[\pi,u_{\ast}]=\endo(J)$ it is highly probable that a distortion map is found. If $\Z[\pi,u_{\ast}]\ne\endo(J)$ it can be a hard problem to prove that some nice candidate is a distortion map, but at least one is able most of the time to find a ``denominator" $m$ such that $m\psi$ lies in the subring $\Z[\pi,u_{\ast}]$; in this case, if $\ell\nmid m$ one can use $m\psi$ as a distortion map on divisors of order $\ell$. Several examples are discussed in \cite{GPRS}.

The aim of this section is to exhibit all supersingular curves of genus $2$ with many automorphisms, describe their automorphisms, and compute the characteristic polynomial of $\pi$, which is always a necessary ingredient in order to analyze the structure of the ring $\Z[\pi,u_{\ast}]$. 
Recall that a curve $C$ is said to have {\it many automorphisms} if it has some geometric automorphism other than the identity and the hyperelliptic involution; in other words, if $|\aut_{\kb}(C)|>2$. 

Igusa found equations for all geometric curves of genus $2$ with many automorphisms, and he grouped these curves in six families according to the possible structure of the automorphism group \cite{Ig}, \cite{iko}. Cardona and Quer found a faithful and complete system of representatives of all these curves up to $\kb$-isomorphism and they gave conditions to ensure the exact structure of the automorphism group of each concrete model \cite{car1}, \cite{car4}. The following theorem sums up these results.

\begin{theorem}\label{list}
Any curve of genus $2$ with many automorphisms is geometrically isomorphic to one and only one of the curves in these six families:

\begin{center}
\begin{tabular}{|c|c|c|c|}
\hline Equation of $C$& &$\aut'_{\kb}(C)$&$\aut_{\kb}(C)$\\
\hline $y^2=x^6+ax^4+bx^2+1$&$a,\,b$ satisfy (\ref{fam1})&$C_2$&$C_2\times C_2$\\
\hline $y^2=x^5+x^3+ax$&$a\ne0,\,1/4,\,9/100$&$C_2\times C_2$&$D_8$\\
\hline $y^2=x^6+x^3+a$&$p\ne3,\ a\ne0,\,1/4,\,-1/50$&$S_3$&$D_{12}$\\
 $y^2=ax^6+x^4+x^2+1$&$p=3,\ a\ne0$&$S_3$&$D_{12}$\\
\hline $y^2=x^6-1$&$p\ne3,\,5$&$D_{12}$&$2D_{12}$\\
\hline $y^2=x^5-x$&$p\ne5$&$S_4$&$\tilde{S}_4$\\
&$p=5$&$\op{PGL}_2(\mathbb{F}_5)$&$\tilde{S}_5$\\
\hline $y^2=x^5-1$&$p\ne5$&$C_5$&$C_{10}$\\
\hline
\end{tabular}
\end{center}
\end{theorem}
\begin{equation}\label{fam1}
(4c^3-d^2)(c^2-4d+18c-27)(c^2-4d-110c+1125)\ne0,\ c:=ab,\,d:=a^3+b^3.
\end{equation}

Ibukiyama-Katsura-Oort determined, using Theorem \ref{sscrit}, when the last three curves are supersingular \cite[Props. 1.11, 1.12, 1.13]{iko}:
$$
\begin{array}{l}
y^2=x^6-1 \mbox{ is supersingular iff }\  p\equiv-1\md3 \\
y^2=x^5-x \mbox{ is supersingular iff }\  p\equiv 5,7\md8 \\
y^2=x^5-1 \mbox{ is supersingular iff }\  p\equiv 2,3,4\md5 
\end{array}
$$   

It is immediate to check that $y^2=ax^6+x^4+x^2+1$ is never supersingular if $p=3$. One can apply Theorem \ref{sscrit} to the other curves in the first three families to distinguish the supersingular ones.

\begin{theorem}\label{noq}
Suppose $q$ is a square and let $C$ be a supersingular curve belonging to one of the first five families of Theorem \ref{list}. Then there a twist of $C$ with Weil polynomial $(x+\sqrt{q})^4$, and this twist is unique. 
\end{theorem}

\begin{proof}Let $E$ be a supersingular elliptic curve defined over $\fp$.
By \cite[Prop. 1.3]{iko} the Jacobian $J$ of $C$ is geometrically isomorphic to the product of two supersingular elliptic curves, which is in turn isomorphic to $E\times E$ by a well-known theorem of Deligne. The principally polarized surface $(J,\Theta)$ is thus geometrically isomorphic to $(E\times E,\la)$ for some principal polarization $\la$. Since $E$ has all endomorphisms defined over $\mathbb{F}_{p^2}$, $(E\times E,\la)$ is defined over $\mathbb{F}_{p^2}$ and by a classical result of Weil it is $\mathbb{F}_{p^2}$-isomorphic to the canonically polarized Jacobian of a curve $C_0$ defined over $\mathbb{F}_{p^2}$. By Torelli, $C_0$ is a twist of $C$. The Weil polynomial of $C_0$ is $(x\pm \sqrt{q})^4$ because the Frobenius polynomial of $E$ is $x^2+p$. The fact that $C_0$ and $C_0'$ are the unique twists of $C_0$ with Weil polynomial $(x\pm \sqrt{q})^4$ is consequence of Proposition \ref{qsq}.\hfill{$\Box$}  
\end{proof}

\begin{corollary}\label{noq2}
Under the same assumptions:
\begin{enumerate}
\item The Weil polynomial of $C$ is $(x\pm\sqrt{q})^4$ if and only if $W=(1)^6$.
\item If $C$ belongs to one of the first three families of Theorem \ref{list}, then it admits no twist with Weil polynomial $x^4\pm qx^2+q^2$ or $x^4+q^2$.
\item If any of the curves $y^2=x^5+x^3+ax$, $y^2=x^6+x^3+a$ is supersingular then $a\in \mathbb{F}_{p^2}$.
\end{enumerate}
\end{corollary}

\begin{proof}
(1) By Table \ref{t4}, the set $W_0$ of Weierstrass points of $C_0$ has $G_k$-structure $W_0=(1)^6$ and Lemma \ref{orbit} shows that for all automorphisms $v\ne1,\iota$ one has $W_v\ne(1)^6$; thus, only the twists $C_0$ and $C_0'$ have $W=(1)^6$.

(2) The geometric automorphisms $v$ of $C_0$ satisfy neither $v^6=1,\,v^2\ne1,\,v^3\ne1,\iota$, nor $v^6=\iota$, nor $v^4=\iota$; thus, by Proposition \ref{qsq} the  Weil polynomial of a twist of $C_0$ is neither $x^4\pm qx^2+q^2$ nor $x^4+q^2$.

(3) The Igusa invariants of $C_0$ take values in $\mathbb{F}_{p^2}$ and $a$ can be expressed in terms of these invariants \cite{car2}.\hfill{$\Box$} 
\end{proof}

In a series of papers Cardona and Quer studied the possible structures of the pointed sets $H^1(G_k,\aut_{\kb}(C))$ and found representatives $v\in\aut_{\kb}(C)$ (identified to $1$-cocycles of $H^1(G_k,\aut_{\kb}(C))$) of the twists of all curves with many automorphisms \cite{car1}, \cite{car2}, \cite{car3}, \cite{car4}. 
In the next subsections we compute the zeta function and the number of $k$-automorphisms of these curves when they are supersingular. A general strategy  that works in most of the cases is to apply the techniques of Sect. \ref{sec2} to find a twist of $C$ with Weil polynomial $(x\pm\sqrt{q})^4$ (for $q$ square) or $(x^2\pm q)^2$ (for $q$ nonsquare) and apply then Propositions \ref{qsq}, \ref{qnsq} to obtain the zeta function of all other twists of $C$.
The results are displayed in the Appendix in the form of Tables, where we exhibit moreover an equation of each curve.

\subsection{Twists of the curve $C\colon y^2=x^5-1$, for $p\not\equiv0,\,1\md5$}
We have $\phi(W)=\{\infty\}\cup \mu_5$ and $\aut'_{\kb}(C)\simeq\mu_5$. The zeta function of $C$ can be computed from Tables 3,4 and Lemma \ref{modauto} applied to $C\otimes k_2$. If $q\not\equiv 1\md5$ the only twists are $C$, $C'$. If $q\equiv1\md5$ there are ten twists and their zeta function can be deduced from Proposition \ref{qsq}. Table \ref{t5} summarizes all computations. 

\subsection{Twists of the curve $C\colon y^2=x^5-x$, for $p\equiv5,\,7\md8$}
Now $\phi(W)=\{\infty,\,0,\,\pm1,\,\pm i\}$. If $p=5$ we have $\aut'_{\kb}(C)=\aut(\pr1)$. If $p\ne5$ the group $\aut'_{\kb}(C)$ is isomorphic to $S_4$ and it is generated by the transformations $T(x)=ix$, $S(x)=\frac{x-i}{x+i}$, with relations
$S^3=1=T^4$, $ST^3=TS^2$. For $q$ nonsquare the zeta function of $C$ is determined by Table \ref{t3}; since the curve is defined over $\mathbb{F}_p$ we obtain the zeta function of $C$ over $k$ by scalar extension.

In all cases we can apply Propositions \ref{qsq} and \ref{qnsq} to determine the zeta function of the twists of $C$. Tables \ref{t6}, \ref{t7}, \ref{t9} summarize all computations. 

\subsection{Twists of the curve $C\colon y^2=x^6-1$, for $p\equiv-1\md3$, $p\ne5$}
We have $\phi(W)=\mu_6$ and $\aut'_{\kb}(C)=\{\pm x,\,\pm\eta x,\pm \eta^2x,\pm\frac1x,\pm\frac{\eta}x,\pm\frac{\eta^2}x\}$, where $\eta\in\mathbb{F}_{p^2}$ is a primitive third root of unity.

The zeta function of $C$ can be computed from Tables 3,4 and Lemma \ref{modauto} applied to $C$ and $C\otimes k_2$. The zeta function of all twists can be determined by Propositions \ref{qsq}, \ref{qnsq}. Tables \ref{t10}, \ref{t11} summarize all computations.   

\subsection{Twists of the supersingular curve $C\colon y^2=x^6+x^3+a$, for $p>3$} Recall that $a$ is a special value making the curve $C$ supersingular and $a\ne 0,\,1/4,\,-1/50$. We have now 
$$\phi(W)=\{\theta,\,\eta\theta,\,\eta^2\theta,\,\frac A{\theta},\,\eta\frac A{\theta},\,\eta^2\frac A{\theta}\},\quad
\aut'_{\kb}(C)=\{x,\,\eta x,\,\eta^2 x,\,\frac Ax,\,\eta\frac Ax,\,\eta^2\frac Ax\},$$
where $A,\,z,\,\theta\in\kb$ satisfy 
$A^3=a, \ z^2+z+a=0,\ \theta^3=z$. 

The Galois action on $W$ and on $\aut'_{\kb}(C)$ depends on $z$ and $a/z$ being cubes or not in their minimum field of definition $k^*$ or $(k_2)^*$. This is determined by the fact that $a$ is a cube or not.

\begin{lemma}
If $a$ is a cube in $k^*$ then $z,\,a/z$ are both cubes in $k^*$ or in $(k_2)^*$, according to $1-4a\in (k^*)^2$ being a square or not. 

If $a$ is not a cube in $k^*$ then $z,\,a/z$ are both noncubes in $k^*$ or in $(k_2)^*$, according to $1-4a\in (k^*)^2$ being a square or not. 
\end{lemma}

\begin{proof}
Let us check that all situations excluded by the statement lead to $W=(1)^3(3)$ or $W_v=(1)^3(3)$ for some twist, in contradiction with Tables \ref{t3}, \ref{t4}. 

Suppose $q\equiv-1\md3$. If $1-4a$ is a square then $a,\,z,\,a/z$ are all cubes in $k^*$. If $1-4a$ is not a square then $a$ is a cube and if $z,z^{\sg}$ are not cubes in $k_2$ we have $\theta^{\sg}=\omega (A/\theta)$, with $\omega^3=1,\,\omega\ne1$, and the twist by $v=(\omega^{-1}(A/x),\sqrt{a}y/x^3)$ has $W_v=(1)^3(3)$ by Lemma \ref{orbit}.  

Suppose $q\equiv1\md3$. If $1-4a$ is not a square we have $z^{(q^2-1)/3}=a^{(q-1)/3}$, so that $a$ is a cube in $k^*$ if and only if $z,\,z^{\sg}$  are cubes in $k_2^*$. Suppose now that $1-4a$ is a square. If exactly one of the two elements $z,\,a/z$ is a cube we have $W=(1)^3(3)$; thus $z,\,a/z$ are both cubes or noncubes in $k^*$. In particular, if $a$ is not a cube then $z,\,a/z$ are necessarily both noncubes. Finally, if $a$ is a cube and $z,\,a/z$ are noncubes in $k^*$, Lemma
\ref{orbit} shows that $W_v=(1)^3(3)$ for the twist corresponding to $v=(\eta x,y)$.\hfill{$\Box$}
\end{proof}

For the computation of the zeta functions of the twists it is useful to detect that some of the combinations $a$ square/nonsquare and $1-4a$ square/nonsquare are not possible.

\begin{lemma}Suppose $q\equiv1\md3$.
\begin{enumerate}
\item If $q\equiv -1\md4$ then $1-4a$ is not a square.
\item If $q$ is nonsquare then $a$ is not a square.
\item If $q$ is a square then $a$ and $1-4a$ are both squares.
\end{enumerate}
\end{lemma}

\begin{proof}Let $C_v$ be the twist of $C$ corresponding to $v(x,y)=(\eta x,y)$.

(1) Supose $1-4a$ is a square. If $a$ is a cube we have $W=(1)^6$ and if $a$ is not  a cube we have $W_v=(1)^6$; by Table \ref{t3} we get $(r,s)=(0,-2\sgn{p}q)$ in both cases. On the other hand, Lemmas \ref{modauto} and
\ref{orbit} applied to $C\otimes_kk_2$ show in both cases that $s\equiv1\md3$; thus, $p\equiv1\md4$.  

(2)  Suppose $a$ is a square. If $a$ is a cube (respectively $a$ is not a cube) we have $W=(1)^6$ or $W=(2)^3$ (respectively $W_v=(1)^6$ or $W_v=(2)^3$), according to $1-4a$ being a square or not. In all cases we have $(r,s)=(0,\pm 2q)$ by Table \ref{t3}, and a straightforward application of Lemma \ref{modauto} and (2) of Lemma \ref{orbit} leads to $r\equiv-1\md3$, which is a contradiction.

(3) In all cases in which $a$ or $1-4a$ are nonsquares we get $(r,s)=(0,q)$ either for the curve $C$ or for the curve $C_v$. This contradicts Corollary \ref{noq2}.  \hfill{$\Box$}
\end{proof}

After these results one can apply the general strategy. The results are displayed in Tables \ref{t12}, \ref{t13}, \ref{t14}. 

\subsection{Twists of the supersingular curve $C\colon y^2=x^5+x^3+ax$}
Recall that $a$ is a special value making $C$ supersingular and $a\ne0,\,1/4,\,9/100$. Given $z\in\kb$ satisfying $z^2+z+a=0$ we have 
$\phi(W)=\{0,\,\infty,\,\pm\sqrt{z},\,\pm\sqrt{a/z}\}$,
$$\aut_{\kb}(C)=\left\{(\omega^2\, x,\omega\, y)\tq \omega^4=1\right\}\cup \left\{\left(\frac{w^2}x,\frac {w^3y}{x^3}\right)\tq w^4=a\right\}\enspace.$$ 

\begin{lemma}\label{15}
If $q\equiv 1\md4$  then $a$ and $1-4a$ are both squares
or both nonsquares in $k^*$.  
If $q$ is a square then necessarily $a$ and $1-4a$ are both squares.
\end{lemma}

\begin{proof}
 If $a\not\in (k^*)^2$, $1-4a\in (k^*)^2$, then $W=(1)^4(2)$ and $(r,s)=(0,-2q)$ by Tables \ref{t3},\ref{t4}; this contradicts 
 Lemma \ref{modauto} because $|\aut(C)|=|\f|=4$ and $r\equiv2\md4$. 
 
Suppose now $a\in (k^*)^2$, $1-4a\not \in (k^*)^2$. If $a\in (k^*)^4$ then $W=(1)^2(2)^2$ and $(r,s)=(0,\pm2q)$; this contradicts 
Lemma \ref{modauto} because $|\aut(C)|=8$, $|\f|=6$ if $q\equiv1\md8$ and $|\f|=2$ or $10$ if $q\equiv5\md8$, so that $r\equiv4\md8$ in both cases. If $a\not\in(k^*)^4$ we get a similar contradiction for the curve $C_v$ for $v(x,y)=(-x,iy)$.

If $a$, $1-4a$ are nonsquares, then $W=(1)^2(4)$ and the Weil polynomial of $C$ is $x^4+q^2$ by Tables \ref{t3},\ref{t4}. If $q$ is a square this contradicts Corollary \ref{noq2}. \hfill{$\Box$}
\end{proof}

\begin{lemma}
If $q$ is a square then $a\in (k^*)^4$ if and only if $z\in (k^*)^2$. 
\end{lemma}
  
\begin{proof}
Suppose $a\in (k^*)^4$, $z\not \in (k^*)^2$ and let us look for a contradiction. Consider the $k$-automorphisms $u(x,y)=(-x,iy)$, $v(x,y)=(\frac{w^2}x,\frac{w^3}{x^3}y)$ of $C$, where $w^4=a$. By Lemma \ref{orbit}, $W_u=(1)^6$ and $C_u$ has Weil polynomial $(x\pm\sqrt{q})^4$ by Corollary \ref{noq2}; since $u^2=\iota$, the Weil polynomial of $C$ is $(x^2+q)^2$ by Proposition \ref{qsq} and Lemma \ref{trans}. The quotient $E:=C/v$ is an elliptic curve defined over $k$ and the Frobenius endomorphism $\pi$ of $E$ must satisfy $\pi^2=-q$. Since $q$ is a square, $E$ has four automorphisms and its $j$ invariant is necessarily $j_E=1728$. Now, $E$ has a Weierstrass equation: $Y^2=(X+2w)(X^2+1-2w^2)$, where $X=(x^2+w^2)x^{-1}$, $Y=y(x+w)x^{-2}$ are invariant under the action of $v$. The condition $j_E=1728$ is equivalent to $a=0$ (which was excluded from the beginning) or $a=(9/14)^2$; in this latter case $z$ is a square in  $\mathbb{F}_{p^2}$ and we get a contradiction.

Suppose now $a\not \in (k^*)^4$, $z\in (k^*)^2$. We have $W=(1)^6$ and $C$ has Weil polynomial $(x\pm \sqrt{q})^4$ by Corollary \ref{noq2}. By Proposition \ref{qsq}, the Weil polynomial of $C_u$ is $(x^2+q)^2$. For any choice of $w=\root4\of{a}$, the morphism $f(x,y)=(\frac{x+w}{x-w},\frac{8\sqrt{w^3}}{\sqrt{1+2w^2}}\frac{y}{(x-w)^3})$ sets a $k_2$-isomorphism between $C$ and the model: 
$$C_u\colon y^2=(x^2-1)(x^4+bx^2+1),\qquad b=(12\sqrt{a}-2)/(2\sqrt{a}+1),    
$$of $C_u$. The quotient of this curve by the automorphism $(-x,y)$ is the elliptic curve $E\colon Y^2=(X-1)(X^2+bX+1)$.  Arguing as above, $E$ has $j$-invariant $1728$, and this leads to $a=0$ (excluded from the beginning) or $a=(9/14)^2$, which is a contradiction since $a$ would be a fourth power in $\mathbb{F}_{p^2}$.\hfill{$\Box$}
\end{proof}

After these results one is able to determine the zeta function of all twists of $C$ when $q$ is a square; the results are displayed in Table \ref{t15}. In the cases where the Weil polynomial is $(x-\eps\sqrt{q})^4$, $\eps=\pm1$, the methods of section \ref{sec2} are not sufficient to determine $\eps$; our computation of this sign follows from a study of the $4$-torsion of an elliptic quotient of the corresponding curve.

In order to deal with the case $q$ nonsquare we need to discard more cases. 

\begin{lemma}Suppose $q$ nonsquare. If $q\equiv -1\md4$  then $a$ and $1-4a$ cannot be both nonsquares. 

If $q\equiv 1\md4$ and $a\in (k^*)^2$ then $a\in (k^*)^4$ if and only if $z\not \in (k^*)^2$. 
\end{lemma}

\begin{proof}
If $a$, $1-4a$ are both nonsquares the polynomial $x^4+x^2+a$ is irreducible and the Weil polynomial of $C$ is $x^4+q^2$ by Table \ref{t3}; hence, the Weil polynomial of $C\otimes_kk_2$ is $(x^2+q^2)^2$. If $q\equiv-1\md4$ we have $a\in k^*\subseteq (k_2^*)^4$ and this contradicts Table \ref{t15}.

Suppose $q\equiv 1\md4$ and $a\in (k^*)^2$; by Lemma \ref{15}, $1-4a$ is also a square and $z\in k^*$. If $a\in (k^*)^4$ and $z\in (k^*)^2$
we get $W=(1)^6$, and $(r,s)=(0,-2q)$ by Table \ref{t3}; we get a contradiction because the Jacobian $J$ of $C$ is simple (\cite[Thm. 2.9]{mn}) and $C$ has elliptic quotients over $k$ because the automorphisms $(w^2/x,(w^3y)/x^3)$ are defined over $k$.  If $a\not \in (k^*)^4$ and $z\not \in (k^*)^2$
we get an analogous contradiction for the curve $C_u$ twisted by $u(x,y)=(-x,iy)$.\hfill{$\Box$}
\end{proof}

The results for the case $q$ nonsquare follow now by the usual arguments and they are displayed in Tables \ref{t16}, \ref{t17}.

\subsection{Twists of the supersingular curve $C\colon y^2=x^6+ax^4+bx^2+1$}
Recall that $a,\,b\in k$ are special values satisfying (\ref{fam1}) and making $C$ supersingular; in particular $p>3$.
The curve $C$ has four twists because $\aut_{\kb}(C)=\aut(C)=\{(\pm x,\pm y)\}$ is commutative and has trivial Galois action.
The Jacobian of $C$ is $k$-isogenous to the product $E_1\times E_2$ of the elliptic curves with Weierstrass equations $y^2=x^3+ax^2+bx+1$, $y^2=x^3+bx^2+ax+1$, obtained as the quotient of $C$ by the respective automorphisms $v=(-x,y)$, $\iota v=(-x,-y)$. 
For $q$ nonsquare, these elliptic curves have necessarily Weil polynomial $x^2+q$ and the Weil polynomial of $C$ is $(x^2+q)^2$. 

\begin{lemma}
If  $q$ is a square $C$ has Weil polynomial $(x\pm\sqrt{q})^4$.
\end{lemma}

\begin{proof}
By Theorem \ref{noq} and Proposition \ref{qsq} $C$ has Weil polynomial $(x\pm\sqrt{q})^4$ or $(x^2-q)^2$. In both cases the elliptic curves $E_1$, $E_2$ have Weil polynomial $(x\pm\sqrt{q})^2$ and we claim that they are isogenous. Since $E(k)\simeq (\Z/(1\pm\sqrt{q})\Z)^2$ as an abelian group, our elliptic curves have four rational $2$-torsion points and the polynomial $x^3+ax^2+bx+1$ has three roots $e_1,\,e_2,\,e_3\in k$. Since $e_1e_2e_3=1$, either one or three of these roots are squares. If only one root is a square we have $W=(1)^2(2)^2$, $W_v=(1)^4(2)$ and $C$, $C_v$ have both Weil polynomial $(x^2\pm q)^2$, in contradiction with Theorem \ref{noq}. Hence, the three roots are squares, $W=(1)^6$, and $C$ has Weil polynomial $(x\pm\sqrt{q})^4$ by Corollary \ref{noq2}.\hfill{$\Box$} 
\end{proof}
 
The zeta function of the twists of $C$ is obtained from Propositions \ref{qsq} and \ref{qnsq}. The results are displayed in Table \ref{t18}. For $q$ square the sign of $(x\pm\sqrt{q})^4$ can be determined by analyzing the $4$-torsion of the elliptic curve
$y^2=x^3+ax^2+bx+1$.   

Finally, there are special curves over $k$ whose geometric model $y^2=x^6+ax^4+bx^2+1$ is not defined over $k$ (cf. \cite[Sect.1]{car1}). It is straightforward to apply the techniques of this paper to determine their zeta function too. 

\section{Appendix}
In this appendix we display in several tables the computation of the zeta function of the supersingular curves of genus $2$ with many automorphisms. For each curve $C_v$, we exhibit the number of $k$-automorphisms and the pair of integers $(r,s)$ determining the Weil polynomial $f_{J_v}(x)=x^4+rx^3+sx^2+qrx+q^2$ of $C_v$. In the column labelled ``s.d" we indicate if $C$ is self-dual. For the non-self-dual curves we exhibit only one curve from the pair $C_v$, $C_v'$. 

We denote by $\eta,\,i\in \kb$ a primitive third,\,fourth root of unity. For $n$ a positive integer and $x\in k^*$ we define $$\nu_n(x)=1 \ \mbox{ if } x\in(k^*)^n, \qquad  \nu_n(x)=-1 \ \mbox{ otherwise}\enspace.$$ In all tables the parameters $s,\,t$ take values in $k^*$.

\begin{scriptsize}
\begin{table}
\caption{\scriptsize Twists of the curve $y^2=x^5-1$ for $p\equiv2,3,4\md5$. The sign $\eps=\pm1$ is determined by $\sqrt{q}\equiv \eps\md5$. The last row provides eight inequivalent twists corresponding to the four nontrivial values of $t\in k^*/(k^*)^5$}
\label{t5}
\begin{center}
\begin{tabular}{|c|c|c|c|c|c|}
\hline  $C_v$&$v$&  &$(r,s)$&\mbox{\scriptsize s.d.}&\mbox{\scriptsize $|\aut(C_v)|$}\\
\hline $y^2=x^5-1$&$(x,y)$&$\begin{array}{c}q\equiv\pm2\md5\\q\equiv-1\md5\\q\equiv1\md5\end{array}$&$\begin{array}{c}(0,0)\\(0,2q)\\(-4\eps\sqrt{q},6q)
\end{array}$&no&$\begin{array}{c}2\\2\\10
\end{array}$\\
\hline$y^2=tx^5-1$, \ \mbox{\scriptsize $t\not\in (k^*)^5$}&$(t^{\frac{1-q}5}x,y)$&$q\equiv1\md5$&$(\eps\sqrt{q},q)$&no&$10$\\
\hline
\end{tabular}
\end{center}
\end{table}
\begin{table}
\caption{\scriptsize Twists of the curve $y^2=x^5-x$ when $q\equiv-1\md8$}
\label{t6}
\begin{center}
\begin{tabular}{|c|c|c|c|c|}
\hline  $C_v$&$v$&$(r,s)$&\mbox{\scriptsize s.d.}&\mbox{\scriptsize $|\aut(C_v)|$}\\
\hline $y^2=x^5-x$&$(x,y)$&$(0,2q)$&yes&$8$\\
\hline $y^2=x^5+x$&$(ix,\frac{1+i}{\sqrt2} y)$&$(0,2q)$&yes&$4$\\
\hline$\begin{array}{c}y^2=(x^2+1)(x^2-2tx-1)(x^2+\frac 2tx-1)\\\mbox{\scriptsize $t^2+1\not\in (k^*)^2$}\end{array}$&$(-\frac1x,\frac y{x^3})$&$(0,-2q)$&yes&$24$\\
\hline$\begin{array}{c}y^2=(x^2+1)(x^4-4tx^3-6x^2+4tx+1),\\\mbox{\scriptsize $t^2+1\not\in (k^*)^2$}\end{array}$&$\left(\frac ix,\frac{i-1}{\sqrt{2}x^3}y\right)$&$(0,0)$&yes&$4$\\
\hline$\begin{array}{c}y^2=x^6-(t+3)x^5+5(\frac{2+t-s}2)x^4+5(s-1)x^3\\ +5(\frac{2-t-s}2)x^2+(t-3)x+1\ \mbox{\scriptsize irred.,  $s^2+t^2=-2$}
\end{array}$ &$\left(\frac{x-i}{x+i},\frac{2(1-i)y}{(x+i)^3}\right)$&$(0,q)$&no&$6$\\
\hline
\end{tabular}
\end{center}
\end{table}
\begin{table}
\caption{\scriptsize Twists of the curve $y^2=x^5-x$ when $q\equiv 5\md8$}
\label{t7}
\begin{center}
\begin{tabular}{|c|c|c|c|c|c|}
\hline $C_v$&$v$&$p$&$(r,s)$&\mbox{\scriptsize s.d.}&\mbox{\scriptsize $|\aut(C_v)|$}\\
\hline $y^2=x^5-x$&$(x,y)$&$\begin{array}{c}p>5\\p=5\end{array}$&$(0,-2q)$&yes&$\begin{array}{c}24\\120\end{array}$\\
\hline $y^2=x^5-4x$&$(-x,i y)$&&$(0,2q)$&yes&$8$\\
\hline $y^2=x^5-2x$&$(ix,\frac{1+i}{\sqrt2} y)$&&$(0,0)$&yes&$4$\\
\hline$y^2=(x^2+2)(x^4-12x^2+4)$&$\left(\frac ix,\frac{i-1}{\sqrt{2}x^3}y\right)$&$\begin{array}{c}p>5\\p=5\end{array}$&$(0,2q)$&yes
&$\begin{array}{c}4\\12\end{array}$\\
\hline $\begin{array}{c}y^2=f(t,x)f(\frac{18+(5i-3)t}{(5i+3)-2t},x)\\ \mbox{\scriptsize $f(t,x)=x^3-tx^2+(t-3)x+1$ \  irred.}\end{array}$&$\left(\frac{x-i}{x+i},\frac{2(1-i)y}{(x+i)^3}\right)$&$\begin{array}{c}p>5\\p=5\end{array}$&$(0,q)$& 
$\begin{array}{c}\mbox{no}\\\mbox{yes}\end{array}$&$6$\\
\hline $y^2=x^5-x-t$,\quad \mbox{\scriptsize $\op{tr}_{k/\mathbb{F}_5}(t)=1$}&$(x+1,y)$&$p=5$&$(\sqrt{5q},3q)$&no&$10$\\
\hline $y^2=x^6+tx^5+(1-t)x+2,\quad$\mbox{\scriptsize irred.}&$(\frac 3{x-1},\frac{\sqrt{2}y}{(x+1)^3})$&$p=5$&$(0,-q)$&yes&$6$\\
\hline
\end{tabular}
\end{center}
\end{table}
\begin{table}
\caption{\scriptsize Twists of the curve $y^2=x^5-x$ when $p\equiv 5,\,7\md8$ and $q$ is a square. Here $\eps=(-1/\sqrt{q})$ and $\eps'=(-3/\sqrt{q})$}
\label{t9}
\begin{center}
\begin{tabular}{|c|c|c|c|c|c|}
\hline  $C_v$&$v$&$p$&$(r,s)$&\mbox{\scriptsize s.d.}&\mbox{\scriptsize $|\aut(C_v)|$}\\
\hline $y^2=x^5-x$&$(x,y)$&$\begin{array}{c}p>5\\p=5\end{array}$&$(-4\eps\sqrt{q},6q)$&no&$\begin{array}{c}48\\240\end{array}$\\
\hline $y^2=x^5-t^2x$,\ \mbox{\scriptsize $\quad t\not\in (k^*)^2$}&$(-x,i y)$&&$(0,2q)$&yes&$8$\\
\hline $y^2=x^5-tx$,\ \mbox{\scriptsize $\quad t\not\in (k^*)^2$}&$(ix,\frac{1+i}{\sqrt2} y)$&&$(0,0)$&no&$8$\\
\hline$\begin{array}{c}y^2=(x^2-t)(x^4+6tx^2+t^2), \\ \mbox{\scriptsize $t\not\in (k^*)^2$}\end{array}$
&$\left(\frac ix,\frac{i-1}{\sqrt{2}x^3}y\right)$&$\begin{array}{c}p>5\\p=5\end{array}$&$(0,-2q)$&yes&$\begin{array}{c}4\\12\end{array}$\\
\hline$\begin{array}{c}y^2=(x^3-t)(x^3-(15\sqrt{3}-26)t),\\ \mbox{\scriptsize $t\not\in (k^*)^3$}\end{array}$&$\left(\frac{x-i}{x+i},\frac{2(1-i)y}{(x+i)^3}\right)$&
$\begin{array}{c}p>5\\p=5\end{array}$&$(2\eps'\sqrt{q},3q)$&no&$\begin{array}{c}6\\12\end{array}$\\
\hline $y^2=x^5-x-t$,\quad \mbox{\scriptsize $\op{tr}_{k/\mathbb{F}_5}(t)=1$}&$(x+1,y)$&$p=5$&$(\sqrt{q},q)$&no&$10$\\
\hline $y^2=x^6+tx^5+(1-t)x+2,\quad$\mbox{\scriptsize irred.}&$(\frac 3{x-1},\frac{\sqrt{2}y}{(x+1)^3})$&$p=5$&$(0,q)$&no&$12$\\
\hline
\end{tabular}
\end{center}
\end{table}
\begin{table}
\caption{\scriptsize Twists of the curve $y^2=x^6-1$ when $q\equiv-1\md3$, $p\ne5$. Here $\eps=(-1/p)$}
\label{t10}
\begin{center}
\begin{tabular}{|c|c|c|c|c|}
\hline  $C_v$&$v$&$(r,s)$&\mbox{\scriptsize s.d.}&\mbox{\scriptsize $|\aut(C_v)|$}\\
\hline $y^2=x^6-1$&$(x,y)$&$(0,2q)$&iff $\eps=-1$&$6+2\eps$\\
\hline $y^2=x^6-t$,\quad \mbox{\scriptsize $\ t\not\in (k^*)^2$}&$(-x,-y)$&$(0,2q)$&iff $\eps=1$&$6-2\eps$\\
\hline$y^2=x(x^2-1)(x^2-9)$&$(\frac1x,\frac{iy}{x^3})$&$(0,-2\eps q)$&yes&$12$\\
\hline$\begin{array}{c}y^2=(x^4-2stx^3+(7s+1)x^2+2tsx+1)\cdot\\\cdot(x^2-\frac4tx-1), \ \mbox{\scriptsize $t^2+4\in k^*\setminus(k^*)^2,\ s^{-1}=t^2+3$}
\end{array}$&$(-\frac1x,\frac{iy}{x^3})$&$(0,2\eps q)$&yes&$12$\\
\hline$\begin{array}{c}y^2=x^6+6tx^5+15sx^4+20tsx^3+15s^2x^2+\\+6ts^2x+s^3, \
\mbox{\scriptsize $s=t^2-4\not\in  (k^*)^2$},\\ \mbox{\scriptsize $\op{gcd}(x^{(q+1)/3}-1,x^2-tx+1)=1$}\end{array}$&$(\frac{\eta}x,\frac{iy}{x^3})$&$(0,\eps q)$&yes&$6$\\
\hline$\begin{array}{c}y^2=x^6+6x^5+15sx^4+20sx^3+15s^2x^2+\\+6s^2x+s^3, \
\mbox{\scriptsize $s=t^2/(t^2+4)\not\in  (k^*)^2$}, \\ \mbox{\scriptsize $\op{gcd}(x^{(q+1)/3}+1,x^2-tx-1)=1$}\end{array}$&$(-\frac{\eta}x,\frac{iy}{x^3})$&$(0,-\eps q)$&yes&$6$\\
\hline
\end{tabular}
\end{center}
\end{table}
\begin{table}
\caption{\scriptsize Twists of the curve $y^2=x^6-1$ when $p\equiv-1\md3$, $p\ne5$ and $q$ is a square. Here $\eps=(-3/\sqrt{q})$}
\label{t11}
\begin{center}
\begin{tabular}{|c|c|c|c|c|}
\hline  $C_v$&$v$&$(r,s)$&\mbox{\scriptsize s.d.}&\mbox{\scriptsize $|\aut(C_v)|$}\\
\hline $y^2=x^6-1$&$(x,y)$&$(-4\eps\sqrt{q},6q)$&no&$24$\\
\hline $y^2=x^6-t^3$, \ \mbox{\scriptsize $\quad t\not\in (k^*)^2$}&$(-x,y)$&$(0,-2q)$&yes&$12$\\
\hline$y^2=x^6-t^2$, \ \mbox{\scriptsize $\quad t\not\in(k^*)^3$}&$(\eta x,y)$&$(2\eps\sqrt{q},3q)$&no&$12$\\
\hline$y^2=x^6-t$, \ \mbox{\scriptsize $\quad t\not\in ((k^*)^2\cup(k^*)^3)$}&$(-\eta x,-y)$&$(0,q)$&no&$12$\\
\hline$y^2=x(x^2+3t)(x^2+\frac t3)$, \ \mbox{\scriptsize $\quad t\not\in(k^*)^2$}&$(\frac 1x,\frac{iy}{x^3})$&$(0,2q)$&yes&$4$\\
\hline$y^2=x^6+15tx^4+15t^2x^2+t^3$,\ \mbox{\scriptsize $\quad t\not\in(k^*)^2$}&$(-\frac 1x,\frac{iy}{x^3})$&$(0,-2q)$&yes&$4$\\
\hline
\end{tabular}
\end{center}
\end{table}
\begin{table}
\caption{\scriptsize Twists of the supersingular curve $y^2=x^6+x^3+a$, $a\ne0,1/4,-1/50$, when $q\equiv-1\md3$. Here $\eps=\nu_2(a)$ and $A$ is the cubic root of $a$ in $k$}
\label{t12}
\begin{center}
\begin{tabular}{|c|c|c|c|c|}
\hline  $C_v$&$v$  &$(r,s)$&\mbox{\scriptsize s.d.}&\mbox{\scriptsize $|\aut(C_v)|$}\\
\hline $y^2=x^6+x^3+a$&$(x,y)$&$(0,2q)$&iff $\eps=-1$&$3+\eps$\\
\hline$\begin{array}{c}y^2=\theta^{-3}(x-\theta)^6-g(x)^3+a\theta^{3}(x-\theta^{\sg})^6\\
\mbox{\scriptsize $g(x)$ min. polyn. of $\theta\in k_2\setminus k$, $\op{N}_{k_2/k}(\theta)=A^{-1}$}\end{array}$&$(\frac Ax,\frac{\sqrt a}{x^3}y)$&$(0,2\eps q)$&iff $\eps=-1$&$9+3\eps$\\
\hline$\begin{array}{c}y^2=\theta(x-\eta)^6-g(x)^3+a\theta^{-1}(x-\eta^2)^6\\
\mbox{\scriptsize$g(x)=x^2+x+1$, $\theta\in k_2\setminus(k_2^*)^3,\ \op{N}_{k_2/k}(\theta)=a$}\end{array}$&$(\eta\frac Ax,\frac{\sqrt a}{x^3}y)$&$(0,-\eps q)$&no&$6$\\
\hline
\end{tabular}
\end{center}
\end{table}
\begin{table}
\caption{\scriptsize Twists of the supersingular curve $y^2=x^6+x^3+a$, $a\ne0,1/4,-1/50$, when $q\equiv1\md3$ and $q$ is nonsquare. Here $A$ is a cubic root of $a$ in $k$ and $n=3$, if $a\in (k^*)^3$, whereas $A=a$, $n=1$, if $a\not\in (k^*)^3$}
\label{t13}
\begin{center}
\begin{tabular}{|c|c|c|c|c|c|}
\hline  $C_v$&$v$&\mbox{\scriptsize$\nu_3(a)$}&$(r,s)$&\mbox{\scriptsize s.d.}&\mbox{\scriptsize $|\aut(C_v)|$}\\
\hline $y^2=x^6+x^3+a$&$(x,y)$&$\begin{array}{c}\ \,\,1\\\!-1 \end{array}$&$\begin{array}{c}(0,-2q)\\(0,q) \end{array}$&
$\begin{array}{c}\mbox{yes}\\\mbox{no} \end{array}$&$6$\\
\hline $\begin{array}{c}y^2=x^6+tx^3+t^2a,\ \,t\not\in (k^*)^3\\
y^2=x^6+ax^3+a^3\end{array}$&$(t^{\frac{q-1}3} x,y)$&$\begin{array}{c}\ \,\,1\\\!-1 \end{array}$
&$\begin{array}{c}(0,q)\\(0,-2q) \end{array}$&$\begin{array}{c}\mbox{no}\\\mbox{yes} \end{array}$&$6$\\
\hline$\begin{array}{c}y^2=\theta^{-n}(x-\theta)^6-g(x)^3+a\theta^n(x-\theta^{\sg})^6\\
\mbox{\scriptsize $g(x)$ min. polyn. of $\theta\in k_2\setminus k$, $\op{N}_{k_2/k}(\theta)=A^{-1}$}\end{array}$&$(\frac {\root3\of{a}}x,\frac{\sqrt a}{x^3}y)$&&$(0,2q)$&yes&$2$\\
\hline
\end{tabular}
\end{center}
\end{table}
\begin{table}
\caption{\scriptsize Twists of the supersingular curve $y^2=x^6+x^3+a$, $a\ne0,1/4,-1/50$, when $q$ is a square. Here $\eps=(-3/\sqrt{q})$. Also, $A$ is a cubic root of $a$ in $k$ and $n=3$, if $a\in (k^*)^3$, whereas $A=a$, $n=1$, if $a\not\in (k^*)^3$}
\label{t14}
\begin{center}
\begin{tabular}{|c|c|c|c|c|c|}
\hline  $C_v$&$v$&\mbox{\scriptsize$\nu_3(a)$}&$(r,s)$&\mbox{\scriptsize s.d.}&\mbox{\scriptsize $|\aut(C_v)|$}\\
\hline $y^2=x^6+x^3+a$&$(x,y)$&$\begin{array}{c}\ \,\,1\\\!-1 \end{array}$&$\begin{array}{c}(-4\eps\sqrt{q},6q)\\ (2\eps\sqrt{q},3q)\end{array}$&no&$\begin{array}{c}12\\6\end{array}$\\
\hline $\begin{array}{c}y^2=x^6+tx^3+t^2a,\ \,t\not\in (k^*)^3\\
y^2=x^6+ax^3+a^3\end{array}$&$(t^{\frac{q-1}3} x,y)$&$\begin{array}{c}\ \,\,1\\\!-1 \end{array}$&$\begin{array}{c}(2\eps\sqrt{q},3q)\\(-4\eps\sqrt{q},6q) \end{array}$&no&$\begin{array}{c}6\\12\end{array}$\\
\hline$\begin{array}{c}y^2=\theta^{-n}(x-\theta)^6-g(x)^3+a\theta^n(x-\theta^{\sg})^6\\
\mbox{\scriptsize $g(x)$ min. polyn. of $\theta\in k_2\setminus k$, $\op{N}_{k_2/k}(\theta)=A^{-1}$}\end{array}$&$(\frac {\root3\of{a}}x,\frac{\sqrt a}{x^3}y)$&&$(0,-2q)$&no&$4$\\
\hline
\end{tabular}
\end{center}
\end{table}
\begin{table}
\caption{\scriptsize Twists of the supersingular curve $y^2=x^5+x^3+ax$, $a\ne0,1/4,9/100$, when $q$ is a square. The last row provides two inequivalent twists according to the two values of $\sqrt{a}$. Here $\eps=-(-1/\sqrt{q})\nu_4(z)$ and $\eps'=-(-1/\sqrt{q})\nu_4(tz)$, where $z^2+z+a=0$}
\label{t15}
\begin{center}
\begin{tabular}{|c|c|c|c|c|c|}
\hline  $C_v$&$v$&\mbox{\scriptsize$\nu_4(a)$}&$(r,s)$&\mbox{\scriptsize s.d.}&\mbox{\scriptsize $|\aut(C_v)|$}\\
\hline $y^2=x^5+x^3+ax$&$(x,y)$&$\begin{array}{c}\ \,\,1\\\!-1\end{array}$&$\begin{array}{c}(4\eps\sqrt{q},6q)\\(0,2q)\end{array}$&
$\begin{array}{c}\mbox{no}\\\mbox{yes}\end{array}$&$\begin{array}{c}8\\4\end{array}$\\
\hline $y^2=x^5+tx^3+at^2x,\quad \ t\not\in(k^*)^2$&$(-x,t^{\frac{q-1}4}y)$&$\begin{array}{c}\ \,\,1\\\!-1\end{array}$&$\begin{array}{c}(0,2q)\\(4\eps'\sqrt{q},6q)\end{array}$&
$\begin{array}{c}\mbox{yes}\\\mbox{no}\end{array}$&$\begin{array}{c}4\\8\end{array}$\\
\hline $\begin{array}{c}y^2=g(x)\left(\theta^2(x-\theta^{\sg})^4+g(x)^2+\qquad\right.\\\qquad\left.+a\theta^{-2}(x-\theta)^4\right), \  
\mbox{\scriptsize $\ \op{N}_{k_2/k}(\theta)=\sqrt{a}$}\\\mbox{\scriptsize $g(x)$ min. polyn. of $\theta\in k_2\setminus k$}\end{array} $&$(\frac{\sqrt{a}}x,\frac{\root4\of{a^3}}{x^3}y)$&&$(0,-2q)$
&yes&$4$\\
\hline
\end{tabular}
\end{center}
\end{table}
\begin{table}
\caption{\scriptsize Twists of the supersingular curve $y^2=x^5+x^3+ax$, $a\ne0,1/4,9/100$, when $q$ is nonsquare and $a\not\in (k^*)^2$}
\label{t16}
\begin{center}
\begin{tabular}{|c|c|c|c|c|c|}
\hline  $C_v$&$v$&$(-1/p)$&$(r,s)$&\mbox{\scriptsize s.d.}&\mbox{\scriptsize $|\aut(C_v)|$}\\
\hline $y^2=x^5+x^3+ax$&$(x,y)$&$\begin{array}{c}\ \,\,1\\\!-1\end{array}$&$\begin{array}{c}(0,0)\\(0,2q)\end{array}$&
$\begin{array}{c}\mbox{no}\\\mbox{yes}\end{array}$&$\begin{array}{c}4\\2\end{array}$\\
\hline $\begin{array}{c}y^2=(x^2-a)\left(\theta(x-\sqrt{a})^4+(x^2-a)^2+\right.\\
\left.+a\theta^{-1}(x+\sqrt{a})^4\right),\ \mbox{\scriptsize $\theta\in k_2,\ \op{N}_{k_2/k}(\theta)=a$}\end{array}$&$(\frac{\sqrt{a}}x,\frac{\root4\of{a^3}}{x^3}y)$&$\begin{array}{c}\ \,\,1\\\!-1\end{array}$&$\begin{array}{c}(0,2q)\\(0,0)\end{array}$
&$\begin{array}{c}\mbox{yes}\\\mbox{no}\end{array}$&$\begin{array}{c}2\\4\end{array}$\\
\hline
\end{tabular}
\end{center}
\end{table}
\begin{table}
\caption{\scriptsize Twists of the supersingular curve $y^2=x^5+x^3+ax$, $a\ne0,1/4,9/100$, when $q$ is nonsquare and $a\in (k^*)^2$. Here $\eps=(-1/p)$. If $p\equiv-1\md4$ we assume that $\sqrt a$ belongs to $(k^*)^2$}
\label{t17}
\begin{center}
\begin{tabular}{|c|c|c|c|c|c|}
\hline  $C_v$&$v$&\mbox{\scriptsize$\nu_4(a)$}&$(r,s)$&\mbox{\scriptsize s.d.}&\mbox{\scriptsize $|\aut(C_v)|$}\\
\hline $y^2=x^5+x^3+ax$&$(x,y)$&$\begin{array}{c}\ \,\,1\\\!-1\end{array}$&$\begin{array}{c}(0,2q)\\(0,-2q)\end{array}$&
$\begin{array}{c}\mbox{iff }\eps=-1\\\mbox{yes}\end{array}$&$\begin{array}{c}6+2\eps\\4\end{array}$\\
\hline $y^2=x^5+tx^3+at^2x,\quad \ t\not\in (k^*)^2$&$(-x,t^{\frac{q-1}4}y)$&$\begin{array}{c}\ \,\,1\\\!-1\end{array}$&$\begin{array}{c}(0,-2\eps q)\\(0,2q)\end{array}$&
$\begin{array}{c}\mbox{yes}\\\mbox{no}\end{array}$&$\begin{array}{c}4\\8\end{array}$\\
\hline $\begin{array}{c}y^2=g(x)\left(\theta^2(x-\theta^{\sg})^4+g(x)^2+\right.\\\left.+a\theta^{-2}(x-\theta)^4\right),\  
\mbox{\scriptsize $\ \op{N}_{k_2/k}(\theta)=\sqrt{a}$}\\
\mbox{\scriptsize $g(x)$ min. polyn. of $\theta\in k_2\setminus k$}\end{array} $
&$(\frac{\sqrt{a}}x,\frac{\root4\of{a^3}}{x^3}y)$&&$(0,2q)$
&$\mbox{iff }\eps=1$&$6-2\eps$\\
\hline$\begin{array}{c}y^2=g(x)\left(\theta^2(x-\theta^{\sg})^4+g(x)^2+\right.\\\left.+a\theta^{-2}(x-\theta)^4\right), \ 
\mbox{\scriptsize $\ \op{N}_{k_2/k}(\theta)=-\sqrt{a}$}\\\mbox{\scriptsize $g(x)$ min. polyn. of $\theta\in k_2\setminus k$}\end{array} $
&$(\frac{-\sqrt{a}}x,\frac{i\root4\of{a^3}}{x^3}y)$&&$(0,2\eps q)$
&yes&$4$\\
\hline
\end{tabular}
\end{center}
\end{table}
\begin{table}
\caption{\scriptsize Twists of the supersingular curve $y^2=x^6+ax^4+bx^2+1$ satisfying (\ref{fam1})}
\label{t18}
\begin{center}
\begin{tabular}{|c|c|c|c|c|c|}
\hline  $C_v$&$v$&&$(r,s)$&\mbox{\scriptsize s.d.}&\mbox{\scriptsize $|\aut(C_v)|$}\\
\hline $y^2=x^6+ax^4+bx^2+1$&$(x,y)$&$\begin{array}{c}q \mbox{ nonsq.}\\q \mbox{ square}\end{array}$&$\begin{array}{c}(0,2q)\\(\pm 4\sqrt{q},6q)\end{array}$&
no&$4$\\
\hline $\begin{array}{c}y^2=x^6+atx^4+bt^2x^2+t^3\\ \mbox{\scriptsize $t\not\in (k^*)^2$}\end{array}$&$(-x,-y)$&$\begin{array}{c}q \mbox{ nonsq.}\\q \mbox{ square}\end{array}$&$\begin{array}{c}(0,2q)\\(0,-2q)\end{array}$&
no&$4$\\
\hline
\end{tabular}
\end{center}
\end{table}
\end{scriptsize}\bigskip

\noindent{\bf Conclusion. }We show that the zeta function of a supersingular curve of genus two is almost determined by the Galois structure of a finite set easy to describe in terms of a defining equation. For curves with many automorphisms this result is refined to obtain a direct (non-algoritmic) computation of the zeta function in all cases. As an application one gets a direct computation of the cryptographic exponent      
of the Jacobian of these curves. Also, the computation of the zeta function is necessary to determine the structure of the endomorphism ring of the Jacobian and to compute distortion maps for the Weil and Tate pairings.   \bigskip

\noindent{\bf Acknowledgement.} It is a pleasure to thank Christophe Ritzenthaler for his help in finding some of the equations of the twisted curves.

\end{document}